\documentclass{article}
\usepackage{graphicx}
\usepackage{amssymb,amsmath}
\usepackage[english]{babel}

\def\bcm{_{^\bullet}}

\def\timehm{\count31=\time \count32=\count31 \divide\count31 by 60
\number\count31 \multiply\count31 by 60 \advance\count32 by
-\count31 :\ifnum\count32<10 0\fi \number\count32}

\def\ep{\varepsilon}

\def\ap{\alpha}

\def\ga{\gamma}
\def\be{\beta}

\newtheorem{thm}{Theorem}

\newtheorem{lem}[thm]{Lemma}

\newtheorem{pro}{Proposici\'on}[section]

\def\mod{\mathop{\rm mod}\nolimits}

\def\Q{\mathbb{Q}}

\newcommand{\Gal}{{\rm Gal}}
\newcommand{\Frob}{{\rm Frob }}
\newcommand{\trace}{{\rm trace}}

\newcommand{\GL}{{\rm GL}}

\newcommand{\cond}{{\rm cond}}

\begin{document}

\title{{\bf Solving Fermat-type equations  via
modular $\mathbb Q$-curves  over polyquadratic fields }}
\author{Luis Dieulefait\thanks{ldieulefait@ub.edu; this research was supported by project MTM2006-04895 and by a Jos\'{e} Castillejo grant, both from MECD, Spain.}\,\, and Jorge Jim\'enez
Urroz\thanks{jjimenez@ma4.upc.edu}
\\
 }
\date{\empty}

 \maketitle

\vskip -20mm
%

\begin{abstract}
We solve the diophantine equations $x^4 + d y^2 = z^p$ for $d=2$
and $d=3$ and any prime $p>349$ and $p>131$ respectively. The method consists in
generalizing the ideas applied by Frey, Ribet and Wiles in the
solution of Fermat's Last Theorem, and by Ellenberg in the
solution of the equation $x^4 + y^2 = z^p$, and we use
$\Q$-curves, modular forms and inner twists. In principle our
method can  be applied to solve this type of equations for other
values of $d$.
\end{abstract}

\section{Introduction}

One of the consequences of the  so called abc conjecture is that,
for any integers $\ap,\be,\ga$, the generalized Fermat equation,
$\ap x^r+\be y^s=\ga z^t$, has only finitely many solutions in
integers, $x, y, z, r, s, t$, satisfying $(x,y,z)=1$ and
$1/r+1/s+1/t< 1$, excepting $1^r+2^3=3^2$. One piece of evidence towards
this fact is the proof of Fermat's Last Theorem, or the case when
$r, s, t$ such that $1/r+1/s+1/t< 1$ are fixed, in which the
existence of only finitely many primitive solutions, i.e. with
$(x,y,z)=1$, is already proved by Darmon and Granville in
\cite{darg}. Hence, in general, it is not strange to expect that a
given diophantine equation as the one in the title,
\begin{equation}\label{dioph}
x^4+dy^2=z^p,
\end{equation}
will not have any non-trivial primitive solution. Indeed, the
main result of this paper is the following theorem related to
cases $d=2$ and $d=3$.
\begin{thm} The diophantine equation $x^4+dy^2=z^p$, does not have
non-trivial primitive solutions $(x,y,z,d,p)$ with $(x,y,z)=1$,
$d=2$ and $p > 349$ prime, or $d=3$ and $p > 131$ prime.
\end{thm}
Since the proof of Fermat's Last Theorem, analogous  results have been
proved to be true for many equations  similar to (\ref{dioph}),
(see \cite{d}, \cite{dm} or \cite{rib4}), using in all these cases
the same kind of  technique, which very briefly is as \-follows:
if a primitive solution to the diophantine equation would exist,
one could construct an
 elliptic curve, associated to that solution, which would have too many good
properties which can not be verified all together and, hence, lead
to a contradiction. In particular, the mod $p$ Galois representation attached to this curve will be unramified at most of the primes of bad reduction of the curve, where $p$ is the prime appearing as an exponent in (1). However, completing the previous outline in a
concrete example is very far from obvious. One of the main
obstacles that appears is that, to prove non-existence of certain
solution $(A,B,C)$, one needs to prove that the associated curve,
 if it exists, should be modular. Then,  a careful search
among the modular forms in the appropriate spaces would give us the
desired contradiction. When the elliptic curve is defined over the
rationals, modularity is a consequence of the papers by Breuil,
Conrad, Diamond, Taylor and Wiles \cite{w}, \cite{tw} and
\cite{bcdt}, proving the Shimura-Taniyama conjecture. However, in
general, it is not possible to associate a rational elliptic curve
to the problem, or it is even more convenient to attach a curve
defined over a number field. This is for example the case in the
papers \cite{d}, \cite{e1} or \cite{d2}. The natural objects to
deal with, when considering a  general number field, $K$, are
$\mathbb Q$-curves, i.e. elliptic curves defined over $K$ that are
isogenous to all their Galois conjugates by Gal$(K/\mathbb Q)$.
Indeed, in \cite{rib3}  Ribet gives a characterization of $\mathbb
Q$-curves as the elliptic curves that are quotients of abelian
varieties of $\GL_2$-type, and these varieties are, according to
the generalized Shimura-Taniyama conjecture, precisely those
constructed by Shimura associated to a newform of some congruence
subgroup $\Gamma_1(N)$. Hence, at least under this hypothesis, one
could find the $L$-series of a $\mathbb Q$-curve as a product of
$L$-series of newforms of this subgroup. In this case we will say
that the curve is modular. In our case, given a primitive solution
to equation (\ref{dioph}) we associate the curve
\begin{equation}\label{ell}
E^{}_{\{A,B,C\}}:=y^2=x^3+4Ax^2+2(A^2+rB)x,
\end{equation}
where $r^2=-d$. The $2$-isogeny in  page $302$ of \cite{sil} makes
the curve $E_{\{{A},{B},C\}}$ a $\mathbb Q$-curve of degree $2$
completely defined over the field $K=\mathbb
Q(\sqrt{-d},\sqrt{-2})$. It is also worthwhile to note  that this
curve, with discriminant $\Delta=512(A^2+rB)C^p$, and
$j$-invariant $j=\frac{64(5A^2-3rB)^3}{C^p(A^2+rB)}$, has
semistable reduction at any odd prime of bad reduction, and any
such prime is in fact a divisor of  $C$.

\

Modularity of $\mathbb Q$-curves has been proved by Ellenberg and
Skinner, \cite{es}, under a local condition at $3$.  We will use
this result, together with  certain generalizations of results of
Manoharmayum on Serre's conjecture in characteristic $7$, to check
that the $\mathbb{Q}$-curves in (\ref{ell}) are indeed modular,
making the first step towards the proof of the theorem. How to get
a contradiction from this fact depends heavily on the concrete
example, and we will describe it, in each of our cases, in the
following sections.

\

Although the strategy could apply to other values of $d$, in
practice, it gets more and more complicated: observe that in
particular our $\mathbb Q$-curve will be defined over a field
depending on $d$, and the method will involve computations on the
Fourier coefficients of modular forms of level growing with $d$.
Also, since the method requires that in some spaces of modular forms none of them has certain inner twists, it may be the case that for some values of $d$ the method can not be applied.\\
On the other hand, the case $d=2$ is somehow easier and will be
treated in a slightly different way. Hence, as we have already
mentioned,  in this paper we will restrict to two particular
choices, $d=2$ and $d=3$, dedicating to each of them a separate
section.

\

{\bf Acknowledgment:} We would like to thank Jordi Quer for his carefull reading of the
paper, in particular the section related to Theorem \ref{quer}. Also, we thank the referee for
his many comments and suggestions, done with care, which we found very helpful.

\

\section{The diophantine equation $x^4+2y^2=z^p.$}
The main result of this section is the following theorem.
\begin{thm} The diophantine equation $x^4+2y^2=z^p$ does not have non-trivial primitive solutions for
$p>349$.
 \end{thm}
 In this case, given a primitive solution $(A,B,C)$ of the
 equation in the title, the associated elliptic curve
$E^{}_{\{A,B,C\}}:=y^2=x^3+4Ax^2+2(A^2+rB)x$, for $r\in\mathbb
Q(\sqrt{-2})$  a squareroot of $-2$,  is a modular $\mathbb
Q$-curve completely defined over the imaginary quadratic field $K
= \Q(\sqrt{-2})$. Observe that, since the field of complete
definition is unramified at $3$, and the curve has good or
semistable reduction at primes dividing $3$, modularity of the
curve is a direct application of the results in  \cite{es}. We
proceed to solve the diophantine equation by following the idea
introduced by Frey for Fermat's Last Theorem, and already used in
\cite{e1} and \cite{d} for the case of $\Q$-curves over
$\mathbb{Q}(i)$. The curve does have semistable reduction at any
odd prime of bad reduction, and we have to divide in two cases
depending  on  the factorization of $C$.
\vskip5pt

 {\bf Case} (i){\bf : $C$ is not divisible by $3$}.

\vskip5pt

 \noindent
 In this case, from the equation $A^4+2 B^2=C^p $ and the assumption that the solution is primitive, we conclude that $C$ is relatively prime to $6$, thus since $C^p = A^4+2 B^2 \geq 3$, it has an odd prime factor $q>3$.\\
We consider the Weil restriction $V$ of $E_{\{ A,B,C \}}$, a
$\GL_2$-type abelian surface  defined over $\mathbb{Q}$.
Modularity of the $\Q$-curve means modularity of $V$, i.e., of the
two-dimensional Galois representations attached to it. These
Galois representations $\rho_{V, \lambda}$ have field of
coefficients $\mathbb{Q}(\sqrt{2})$ because the $\Q$-curve has
degree $2$ (cf. \cite{e1}, \cite{d2}).  Most of the ramification of these Galois representations is at the primes $q$ dividing $C$, where the curve has semistable reduction, (and therefore,
 ramification is unipotent at any such $q$ provided that it is unramified in $K$, i.e., that $q>2$), and the idea is to consider a prime
 $\pi \mid p$ in $\mathbb{Q}(\sqrt{2})$ and the residual mod $\pi$ representation $\bar{\rho}_{V, \pi}$: since each odd prime appears with
 multiplicity divisible by $p$ in the discriminant of the $\Q$-curve, we know that the restriction to the Galois group of $K$ of
 $\bar{\rho}_{V, \pi }$, and therefore also $\bar{\rho}_{V, \pi }$ itself, will be unramified at any prime $q>2$ (this is basically the same trick used by Frey
 in connection to Fermat's Last Theorem). \\  Therefore, the conductor of $\bar{\rho}_{V, \pi }$ is a power
of $2$, and we want to compute it. First, using Tate's algorithm,
(cf. \cite{sil}), we find the conductor of $E_{\{ A,B,C \}}$ at
the prime $\check{2}\in K$ dividing $2$  to be either
$\check{2}^{10}$ or $\check{2}^{12}$, depending on the parity of
$B$, (observe that $A$ is always odd because the solution $A,B,C$
is primitive). We then apply the formula in \cite{mi} relating the
conductor of an elliptic curve to the conductor of its Weil
restriction to obtain that the $2$-part of the conductor of $V$ is
either $2^{16}$ or $2^{18}$. Since $V$ is an abelian surface of
$\GL_2$-type, the conductor of $\rho_{V, \pi}$ is the square root
of the conductor of $V$, thus has $2$-part equal to $2^8$ or
$2^9$. We conclude that the conductor of $\bar{\rho}_{V, \pi }$ is
exactly $2^8$ or $2^9$ (as it is well-known that
in presence of wild ramification the conductor does not decrease when reducing mod $p$).\\
At this point, we have almost all ingredients, (as in the proof of
Fermat's last theorem),  to say that $\bar{\rho}_{V, \pi }$
corresponds to a modular form of weight $2$ and level $256$ or
$512$. This is  because $V$ is modular, and also because we can
use lowering-the-level results of Ribet, to ensure that we can
take the modular form of minimal level. We are just lacking one
ingredient: for the above claims to make sense, we have to be sure
that $\bar{\rho}_{V, \pi }$ is irreducible. Since we are in case
(i), there is a prime of semistable reduction, $q$, not dividing
$6$,  and  we can apply Ellenberg's generalization to $\Q$-curves
of results of Mazur (cf. \cite{e1}, Proposition 3.2) and conclude
that for $p >13$ the residual representation is in fact
irreducible. \\ From now on, we assume that $p>13$. We have proved
that there is a newform $f$ of weight $2$  and level $256$ or
$512$ such that one of its attached residual characteristic $p$
representations must be irreducible and isomorphic to
$\bar{\rho}_{V, \pi }$. We consider one by one the newforms in
these two spaces, $22$ in total by \cite{st}, and see that this is
impossible. Whenever the modular form $f$ has complex
multiplication, this is the case treated in \cite{e1} and we will
derive a contradiction from this congruence by using the method
developed  there. Contradiction follows from the fact that for a
newform with complex multiplication, as it is well-known, the
attached $\ell$-adic and residual mod $\ell$ representations have
all ``small image" because they are potentially abelian, more
precisely, the residual image is contained in the normalizer of a
Cartan subgroup. On the other hand, the results of \cite{e1} show
that Galois representations coming from $\Q$-curves defined over a
quadratic field $K$ having bad semistable reduction at some prime
$q\nmid 6$ have large image for every $p > p_0$,  with $p_0$
depending only of the field $K$ and the degree $d$ of the
$\Q$-curve (cf. loc. cit., Propositions 3.4 and 3.6), thus they
can not be equivalent to those attached to a modular form with
Complex Multiplication. In our case the field is $\mathbb Q(\sqrt{-2})$, and $d=2$  and we obtain, by Lemma
\ref{analitischen} below,   $p_0= 349$.\\
 Therefore, it only remains to consider the case of modular forms without complex multiplication. In this case, the contradiction will be derived
 from the following observation. The Galois representation ${\rho}_{V, \pi }$ is attached to $V$, and $V$ is the Weil restriction of a $\Q$-curve.
 As it is well-known (cf. \cite{da2} for a similar situation) this implies
 that ${\rho}_{V, \pi }$ has an ``inner twist", namely, that for any prime $t$
 where ${\rho}_{V, \pi }$ is unramified if we call
 $a_t = \trace \; {\rho}_{V, \pi }(\Frob \; t)$, then it holds:
 $$ a_t^\sigma = \psi(t) a_t \qquad \quad (*)$$
 where $\sigma$ is the involution in $\Gal(\Q(\sqrt{2})/\Q)$ and $\psi$ is the quadratic character corresponding to the extension
 $K/\Q$. In more concrete terms, this means that for any unramified place $t$:\\
 $_{^\bullet}$ $a_t = z \sqrt{2}$, for some rational integer $z$, if $t \equiv 5,7 \mod{8}$, or\\
 $_{^\bullet}$ $a_t \in \mathbb{Z}$ if $t \equiv 1,3 \mod{8}.$\\
 Using this information and the Weil bound   for the coefficients $a_t$ it is enough to derive a contradiction. This method works provided that in the
 spaces that we are considering none of the newforms $f$ satisfy simultaneously the two conditions: $f$ does not have complex multiplication and $f$
 has field of coefficients $\Q(\sqrt{2})$ and an inner twist as in (*). \\
 No newform in the spaces of level $256$ and $512$ satisfies these two conditions and so with a few computations we conclude the proof. In fact we
 only use the eigenvalue $c_3$ of the newforms without complex multiplication in these spaces. Actually, all newforms in the space of level $256$ have complex multiplication, and the newforms of level $512$ without complex multiplication are, up to Galois conjugation, the following five (in the notation of \cite{st}):
 $$ 512A1, 512B1, 512D1, 512E1, 512G1 $$
 the eigenvalue $c_3$ of the first four of them is $c_3= \pm \sqrt{2}$, and the last one has $c_3 = \pm \sqrt{6}$.\\ 
 The required computations are as
 follows:\\
 The mod $p$ congruence between $f$ and   ${\rho}_{V, \pi }$ implies that $a_t \equiv c_t \mod{p}$ for every $t$ unramified.
 Since we are assuming $3 \nmid C$, $V$ has good reduction at $3$, thus we have: $a_3 \equiv c_3 \mod{p} \quad (**)$. \\
 From the inner twist (*) and the Weil bound $| a_3 | \leq 2 \sqrt{3}$ we conclude that the only possible values of $a_3$ are:
 $$ a_3 = 0, \pm 1, \pm 2 \pm 3 $$ 
 Therefore, since $c_3 = \pm \sqrt{2}$ or $\pm \sqrt{6}$ we see that the congruence (**) can not hold whenever $p>13$.
 This, together with the analysis done for the complex multiplication case concludes the proof in case (i) for $p > 349$.\\

 {\bf Case} (ii){\bf : $C$ is divisible by $3$}

\vskip5pt

 \noindent 
 First of all, we have to check that the residual representation $\bar{\rho}_{V, \pi }$  is irreducible. If there is a prime $q>3$
  dividing $C$ this follows as we have already recalled from the results in \cite{e1} for any
  $p>13$.   So here, to prove irreducibility, we only have to consider the case $C= 3^{k}$. In this case the conductor of  ${\rho}_{V, \pi }$ is exactly $3\cdot 2^e$, with $e=8$ or $9$. Now we assume that $p > 17$. The inner twist implies that $a_{17} \in \mathbb{Z}$, and it follows from the Weil bound that $|a_{17}| \leq 8$ (the representation is unramified at $17$ because it only ramifies at $2$, $3$ and $p >17$). Suppose that
  $\bar{\rho}_{V, \pi }$ is reducible. Since the Galois representation is attached to an abelian variety and a weight $2$ modular form, it is well-known that in the residually reducible case one has:
  $$ \bar{\rho}_{V, \pi }^{s.s.} \cong \varepsilon \oplus \varepsilon^{-1} \chi $$
  where $\chi$ is the mod $p$ cyclotomic character and $\varepsilon$ is a character with $\cond(\varepsilon)^2 \mid \cond
   {\rho}_{V, \pi }$.
   In our case, this implies that $\cond(\varepsilon) \mid 16$ and in particular that $\varepsilon(17) = 1$. Thus, evaluating both sides of the congruence at $17$ we obtain: $a_{17} \equiv 1 + 17 \mod{p}$. Using the information we have on $a_{17}$ we see that for any prime $p >23$ this congruence can not hold, thus proving irreducibility for such primes.\\
   Having shown irreducibility, we can handle any case where $3 \mid C$ using raising-the-level results of Ribet. Namely, in such a case the mod $p$
   congruence between the modular Galois representation ${\rho}_{V, \pi }$ and the $p$-adic representation corresponding to a modular form $f$ of level
   $256$ or $512$ implies in particular that one can raise the level of $f$  adding the prime $3$ modulo $p$. As it is well-known, this implies in
   particular that the eigenvalue $c_3$ of $f$ must satisfy:
   $$ c_3 \equiv \pm 4  \mod{p} $$
   and using the values of the eigenvalues $c_3$ of all newforms in the two spaces considered we easily check that this congruence is not satisfied
   by any prime $p > 17$. This gives the desired contradiction in Case (ii).

\section{The equation $x^4+3y^2=z^p$}
This section includes the proof of the following theorem.
\begin{thm}\label{tres} The diophantine equation $x^4+3y^2=z^p$ does not have non-trivial primitive solutions for
$p>131$.
 \end{thm}
Again the proof will rely on the modularity of the associated
elliptic curve $E^{}_{\{A,B,C\}}$, given in (\ref{ell}) with
$r^2=-3$ or, in fact, of its Weil restriction. Indeed, as we
mentioned, $\mathbb Q$-curves are quotients of abelian varieties
of $\GL_2$-type. However, in order to be modular, we want to
consider only special cases of $\Q$-curves, namely such that the
abelian variety $W:=$Res$_{K/\mathbb Q}(E_{A,B,C})$ where $K$ is
the field of complete definition of the curve, is itself of
$\GL_2$-type. In this case  one could find the appropriate level
$N$ of the modular forms by the formula which relates the
conductor of a variety with its restriction of scalars. In the
paper \cite{q2}, Theorem 5.4, Quer gives a sufficient condition
for a given $\mathbb Q$-curve, $C/K$, to have Weil restriction of
$GL_2$-type, in terms of the Schur class of certain cohomology
class. We now include that theorem for convenience of the reader.
For any embedding $s\,:\, K\to\overline{\mathbb Q}$ we can choose an isogeny  $\phi_s\,:\, C\to \,^sC$. 
With this we built the $2$-cocycle for the trivial action of $G_\mathbb Q$ on $\mathbb Q^*$,
\begin{equation}
c_K(s,t)=\phi_s\,^s\phi_t\phi_{st}^{-1},
\end{equation}
and we call $\xi_K(C)\in H^2(K/\mathbb Q,\mathbb Q^*)$ its cohomology class. In page 287 of \cite{q2} it is proved that
$\xi_K(C)$ only depends on $C$. Any two cocycles in $H^2(K/\mathbb Q,\mathbb Q^*)$ can be viewed as taking values into $\overline{\mathbb Q}^*$, considered as a  $\Gal(K/\mathbb Q)$-module with trivial action. We will call the Schur class of 
$c_K$, (or its cohomology class $\xi_K(C)$), to\footnote{In page 304 of \cite{q2} there is a misprint in the definition of Schur class. It should be replaced $H^2(G_\Q,\overline{\mathbb Q}^*)$ by $H^2(G,\overline{\mathbb Q}^*)$.}  the image of $c_K$ into $H^2({K/\mathbb Q},\overline{\mathbb Q}^*)$.

\begin{thm}\label{quer} {\rm (Quer)}
Let $C$ be a $\mathbb Q$-curve without complex multiplication completely defined over a minimal splitting field $K$ such that $\xi_K(C)$ has trivial
Schur class, 
with a quadratic splitting character for $\xi_K(C)$. Then  Res$_{K/\mathbb Q}(C)$ is of $\GL_2$-type.
\end{thm}
In order to apply the theorem in our particular case, we have first to ensure that we start with a non-CM curve. This is the content
of the following lemma.
\begin{lem} The curve $E_{\{{A},{B},C\}}$, with $d=3$,   does not have complex multiplication.
\end{lem}
{\it Proof:}  \, The proof is the same as Proposition $6\text{.}1$ of \cite{py}. Indeed, since
$(3,A)=1$ we can reduce $E_{\{{A},{B},C\}}$ modulo the prime $\frak p\in\mathbb Z[\frac12+\frac12\sqrt{-3}]$ above $3$
to obtain an elliptic curve $E_{\frak p}$ such that $|E_{\frak p}|=2$, and so with  End$(E_{\frak p})\otimes\mathbb Q\simeq \mathbb Q(\sqrt{-2})$.
On the other hand, it is easy to see that $j$ can not be a rational number. The result follows, as in \cite{py},
by noting that  End$(E_{\{{A},{B},C\}})\otimes\mathbb Q\hookrightarrow$End$(E_{\frak p})\otimes\mathbb Q$,
and any elliptic curve over $\bar{\mathbb Q}$ with complex multiplication by a field of class number one has to have
rational $j$-invariant.

\

The next step in order to apply Theorem \ref{quer} is to control the Schur class of the curve $E^{}_{\{A,B,C\}}$. We will
use the notation in \cite{q1}. In this way let $M=\mathbb Q\left(\sqrt{-3}\right)$ and  $N=\mathbb Q\left(\sqrt{-2}\right)$, the fields of definition of
the curve and isogeny respectively, so $MN$ is the field of complete definition as a $\mathbb Q$-curve. From the explicit equation
for the isogeny from $E^{}_{\{A,B,C\}}$ to its Galois conjugate, we can compute the  $2$-cocycle $c_{MN}(s,t)$ and then its sign $c_{MN}^{\pm}$
to obtain
$$
c_{MN}^{\pm}=c_{\ap,\be},\qquad \ap=\{1,-3\}, \quad \be=\{1,-2\}.
$$
Although $\xi_{MN}(E^{}_{\{A,B,C\}})$ does not have trivial Schur
class, in \cite{q1} it is proved the existence of certain
extension of $M$, $K$, and an integer $\gamma\in K$, such that
the twist
\begin{equation}
E^{\gamma}_{\{A,B,C\}}:=\gamma y^2=x^3+4Ax^2+2(A^2+rB)x,
\end{equation}
has $\xi_K(E^{\gamma}_{\{A,B,C\}})$ with trivial Schur class.
Moreover, in the same paper,  it is also given a method to find
both the field $K$ and the twist $\gamma$. In particular $K=LM$, where $L$ is
 a cyclic extension $L/\mathbb Q$ with
Dirichlet character $\ep$ such that Inf$[c_\ep]=(1,3)(-1,-2)$. 
The Dirichlet character for the quadratic field $\mathbb Q(\sqrt 6)$ 
already has this property, and so we can take $L=\Q(\sqrt 6)$. We
now have all the ingredients to apply Theorem 3.1 of \cite{q1} in
order to find $\gamma$. In particular, $L\cap M=\mathbb Q$ and so
$e=1$ and, since $K=MN=ML$, we can take  $a_1=-3$ and $b_1=-2$,
(for precise definitions see \cite{q1} page 195: Decomposition of
K). Now, by choosing $\ap_1=-2+\sqrt{6}$, an element of norm
$N_{K/M}(\ap_1)=-2$, and $\ap_0=\sqrt{-2}+\sqrt{-3}$, it is easy
to find $\beta_{\sigma}$ in (1) of \cite{q1} to obtain, choosing
$x=-\frac 14$ in the same formula, $\gamma=2+\sqrt 6$. Hence, the
curve
$$
E^{\gamma}_{^A,^B,^C}:=\gamma y^2=x^3+4Ax^2+2(A^2+rB)x,
$$
is the quotient of an abelian variety of $\GL_2$-type, concretely, its Weil restriction $W$ to $\Q$. \\
Our next step will be proving  that, in fact, $W$ is modular. Let
us recall that $W$ is an abelian variety of dimension $4$, of
$\GL_2$-type, and such that the corresponding two-dimensional
Galois representations $\rho_{W,\lambda}$ have field of
coefficients $\Q(\sqrt{2}, i)$ and two inner twists (because it
comes  from a $\Q$-curve of degree $2$ over the biquadratic field
$K$, see page 226 of \cite{py} for a similar situation). Hence,
the traces $a_t$ of this compatible family of Galois
representations, at unramified places, satisfy the following
rule:\\
$_{^\bullet}$ $a_t \in \mathbb{Z}$, if $t$ totally splits in $K$.\\
$_{^\bullet}$ $a_t = z \sqrt{2}$, for some integer $z$, if $t$ splits in $L$ and is inert in $M$.\\
$_{^\bullet}$ $a_t = z i$, for some integer $z$, if $t$ splits in $M$ and is inert in $L$.\\
$_{^\bullet}$ $a_t = z \sqrt{-2}$, for some integer $z$, if $t$ is inert in $L$ and $M$.\\
(cf. Theorem 5.4 of \cite{q2})\\

 Since the field $K$ is ramified at $3$, the modularity of $W$ does not follow from an easy application of \cite{es} as in the case 
  $d = 2$. Thus, in order to prove modularity we will
 apply (a slight generalization of) the results of Manoharmayum (cf. \cite{ma} and \cite{jm}) proving Serre's conjecture for the field of $7$ elements.
 Since the abelian variety $W$ has good or semistable reduction at $7$ (because the $\Q$-curve has semistable reduction at every odd prime of bad
 reduction) it is enough  to prove, (cf. \cite{d2003}), that one of the attached residual representations in characteristic $7$ is reducible or modular
 to conclude modularity of $W$. Since modularity is well-known in the case of solvable image, we only have to prove modularity for one of the mod $7$
 representations assuming that it has non-solvable image. Let $\bar{\rho}_7$ be such a residual representation. By solvable base change, it is enough to
 prove modularity of the restriction of   $\bar{\rho}_7$ to the Galois group of $L$, and we know from the inner twists of $W$ that this restriction has
 coefficients in $\mathbb{F}_7$ (because ${2}$ is a square mod $7$). As we will see in the next paragraph (when computing the conductor and inner twists of $W$) this
 restriction is also unramified at $3$. The results in \cite{ma} and \cite{jm} prove modularity in this context, except that they require the condition
 that $3$ is unramified in $L$, which is not satisfied in our case. However, if we inspect the proof of this result we find that using recent modularity
 lifting results of Kisin, (cf. \cite{ki}), and generalizations to ramified number fields by Gee, (cf. \cite{ge}), this restriction is no longer necessary.
 In fact, modularity is proved by showing the existence of an elliptic curve over a solvable extension $X$ of $L$ whose mod $7$ representation agrees
 with the restriction to the Galois group of $X$ of $\bar{\rho}_7$, then it only remains to prove modularity of such elliptic curve. It is shown that
 the curve can be chosen such that the corresponding mod $3$ Galois representation is surjective, and $X$ can be chosen such that $X/L$ is unramified
 at $3$. Then as in Wiles' proof of Fermat's Last Theorem we know that this residual representation is modular, due to results of Langlands and Tunnell.
 The condition $3$ unramified in $L$ (therefore also in $X$) is used by Jarvis and Manoharmayum to apply a modularity lifting result and conclude that
 $E$ is modular, but now using the results of Kisin and Gee we know that even for our field $L$ we can conclude modularity of $E$ and this concludes the
 proof that $W$ is modular.\\
 
 {\bf Remark:} Due to recent works of Khare-Wintenberger, Kisin and the first author, Serre's modularity conjecture has been established in full generality (cf. \cite{kw}, \cite{k} and \cite{ddd}). Serre's conjecture implies (as proved by Ribet in \cite{rib3}) the generalized Shimura-Taniyama conjecture, i.e., the modularity of all abelian varieties over $\Q$ of $\GL_2$ type, thus in particular the modularity of $W$ can be seen as a particular case of this much more general result.

\

To compute the $2$-part and the $3$-part of the conductor of $W$
we first get the conductor of $E^{\gamma}_{^A,^B,^C}$ over $K$,
and then apply the formula in \cite{mi}. Let us stress that in particular the curve $E^{\gamma}_{^A,^B,^C}$ has good reduction at primes of $K$ dividing $3$ (this follows from the fact that the solution $A,B,C$ is primitive, thus $3 \nmid A$), so the $3$-part of the conductor of $W$ comes just from base changing from $K$ to $\Q$ (also, since $K/L$ is unramified at $3$, there is no ramification at $3$ if we just consider the intermediate base change from $K$ to $L$). In this way we obtain
that the $3$-part is $3$ and the $2$-part is $2^e$ with $e=3,5,6$
or $7$. The determinant of the representations $\rho_{W,\lambda}$
is $\varepsilon \chi$ where $\chi$ is the $\ell$-adic cyclotomic
character and $\varepsilon$ the quadratic character corresponding
to the quadratic field $L$. \\
We apply again Frey's trick and, based on the shape of the
discriminant of the curve $E^{\gamma}_{^A,^B,^C}$, we deduce that
for any prime $\pi \mid p$ in the field of coefficients, the
residual representation $\bar{\rho}_{W,\pi}$ loses the
ramification at any prime dividing $C$, which are
primes of semistable reduction of the curve (we also know that $3 \nmid C$ because the solution is primitive). Thus, the conductor
of $\bar{\rho}_{W,\pi}$ is
 $3 \cdot 2^e$ with $e=3,5,6$ or $7$. \\

 Observe that any primitive solution verifying $A^4 + 3 B^2=C^p$ must have $C$ coprime with $3$ (as we already mentioned) and odd,  because if $A$ and $B$ were both odd
 then $A^4 + 3 B^2 \equiv 4 \mod{8}$ which can not be a $p$-th power. Hence $C$ must have some prime divisor
 $q >3$, which will be a prime of semistable reduction of  $E^{\gamma}_{^A,^B,^C}$. Thus, the technical condition needed to apply the results in
 \cite{e1} and conclude that for $p$ larger than a certain bound $\bar{\rho}_{W,\pi}$ is irreducible and its image is not contained in the normalizer
 of a Cartan subgroup hold. In particular, we have irreducibility for $p>13$.
 Thus, for $p>13$, the residual representation $\bar{\rho}_{W,\pi}$ should coincide with the residual representation corresponding to a newform of
 weight $2$, level $3 \cdot 2^e$, and nebentypus $\varepsilon$, for $e=3,5,6$ or $7$.
We also know, by the Shimura correspondence,(see \cite{sh1},
Theorem 7.14, and \cite {rib2}, Proposition 2.3) that the
extension $Q(a_1,a_2,\dots, a_n\cdots)$ generated by the
coefficients of the modular form coincide with the algebra of
endomorphisms of  $W$ which, in this case, is $\mathbb
Q(\sqrt{-1},\sqrt{-2})$ by Theorem 5.4 of \cite{q2}. See also
Theorem 3.1 of \cite{q1} and observe that we already have an
element $\ap_1$ of norm $-2$. Finally we can obtain the Nebentypus
of the form since it coincides, as it was conjectured by Serre in
\cite{se},  with the splitting character $\ep$ of
$\xi_K(E_{\{{A},{B},C\}})$. One can find a proof of this fact in
Theorem 5.12 of \cite{py}. \

\

We eliminate as in the previous example all cases of modular forms with CM (Complex Multiplication) because the results of \cite{e1} ensure large
image for $p$ sufficiently large and this contradicts a congruence with a CM newform. Again Lemma \ref{analitischen} applied for the field
 $\Q(\sqrt{-3})$ gives the result for every $p>131$.\\
The rest of the modular forms (those without CM) in the four spaces that we have to consider are eliminated by using the precise information on the
coefficients $a_t$ due to the existence of two inner twists and, as in the previous section, by comparing a few coefficients and using the Weil bound
we easily show that for any $p>131$ and any newform $f$ with eigenvalues $c_t$ in one of the four spaces considered there is a small prime $t$ such
 that the congruence $a_t \equiv c_t \mod{p}$ is impossible. We have computed the required eigenvalues of these newforms using SAGE (cf. \cite{st2}).\\
 More precisely, we have to consider all newforms of weight $2$, nebentypus $\varepsilon$ and level $24, 96, 192$ and $384$: in the level $24$ and $96$ cases all newforms have CM. In the level $192$ case there are some newforms without CM, and we compute all possible values of $c_5$ for them and we obtain: $c_5 = \pm 2 \sqrt{3}$. In the level $384$ case, there are again some newforms without CM, but they satisfy: $c_7 = \pm 2 \sqrt{-6}$ or $\pm  2 \sqrt{-2}$.\\ 
 Using these eigenvalues we have enough to derive a contradiction:\\
a) Level $192$, nebentypus $\varepsilon$ and weight $2$ without CM: $c_5 =  \pm 2 \sqrt{3}$. From the description of the
inner twists and the Weil bound we know that the coefficient $a_5$ of our representation is of the form $z \sqrt{2}$, for an integer $z$ with
$| z | \leq 3$. Then, for any $p>13$ the congruence $a_5 \equiv c_5 \mod{p}$ does not hold. Observe that $5$ is an unramified prime for
${\rho}_{W,\pi}$ because the curve $E^{\gamma}_{^A,^B,^C}$ can not ramify at $5$ (just observe that $5 \nmid C$ because  $C^p = A^4+3B^2$). \\
b) Level $384$, nebentypus $\varepsilon$ and weight $2$ without CM: $c_7 =  \pm 2 \sqrt{-6}$ or $\pm 2 \sqrt{-2}$. First assume that $7$
is an unramified prime for ${\rho}_{W,\pi}$. Then we know from the inner twists and the Weil bound that
$a_7 = z i$ for an integer $z$ with $|z| \leq 5$. Then, for any $p>17$ the congruence $a_7 \equiv c_7 \mod{p}$ does not hold.\\
Now suppose that ${\rho}_{W,\pi}$ ramifies at $7$ (ramification will be unipotent, $7$ being a prime of bad semistable reduction of
$E^{\gamma}_{^A,^B,^C}$). In this case the congruence between the mod $p$ Galois representation attached to $f$ and
$\bar{\rho}_{W,\pi}$ implies that one can raise the level of $f$ modulo $p$, adding the prime $7$. For this to hold, as it is well-known (condition for level-raising), it is necessary that it
   holds: $c_7 \equiv \pm 8 \mod{p}$. This congruence is not satisfied by any prime $p >11$.\\

\
\section{An explicit bound for the prime $p$.}
In both of the previous sections, to exclude the cases in which
the modular form has complex multiplication, one has to compute a
concrete bound, $p_0$, such that for any given prime $p>p_0$, the
Galois representations coming from $\Q$-curves defined over a
quadratic field $K$ having bad semistable reduction at most places
have large image, and then getting a contradiction.  We compute
this bound exactly as in \cite{e1} for the case considered there.
The difference now is that, in any of the two cases in this paper,
the fields to be considered have character of conductor at most
$8$ which will provide  slightly different bounds than that obtained
in \cite{e1}.

\

Let $\cal F$ be a Petersson-orthogonal basis for
$S_2(\Gamma_0(N))$, $\chi$ a Dirichlet character, and let us write
$$
(a_m,L_\chi)_N=\sum_{f\in \cal F}a_m(f),L(f\otimes\chi,1),
$$
where, $f=\sum_{m\ge o}a_m(f)q^m$. Also,  if $M|N$ we will let
$(a_m,L_\chi)_N^M$ to be the contribution to $(a_m,L_\chi)_N$ of
the forms coming from level $M$. As in \cite{e1}, any prime for
which
the inner product\footnote{The identity (\ref{inner}) is proved in
Lemma 3.12 of \cite{e1}.}
\begin{equation}\label{inner}
(a_1,L_{\chi})_{p^2}^{p-new}=(a_1,L_{\chi})_{p^2}-p(p^2-1)^{-1}(a_1-p^{-1}\chi(p)a_p,L_\chi)_p
\end{equation}
is different from zero, would suffices to guarantee large image of
the residual $\bar\rho_{V,\pi}$ Galois representation. The bound
will follow directly from Theorem 1 of \cite{e2}, and Lemma 3.13 of \cite{e1},  which we now
include for completeness.
\begin{thm}\label{ellen} (Ellenberg) Let $\chi$ be a Dirichlet character of conductor $q$, and  $N$ and integer $N\ge 400$,  $N\nmid q$. Let $\sigma$
be a real number with $q^2/2\pi\le \sigma\le Nq/\log N$.
Then we can write
$$
\left|(a_m,L_\chi)_N\right|=4\pi\chi(m)e^{-2\pi m/\sigma N\log N}- E^{(3)}+E_3-E_2-E_1+(a_m,B(\sigma N\log N)),
$$
where
\begin{eqnarray*}
&&\bcm\, \left|(a_m,B(\sigma N\log N))\right|\le 30 \,(400/399)^{3}e^{2\pi}q^2m^{3/2}N^{-1/2}d(N)N^{-2\pi\sigma/q^2},\\
&&\bcm\, E_1\le\frac{_{16}}{^3}\pi^3m^{3/2}\sigma\log N \,e^{-N/(2\pi m\sigma\log N)},\\
&&\bcm\, E_2\le \frac{_{8}}{^9}\pi^5\zeta^2(7/2)m^{5/2}\sigma^2N^{-3/2}(\log N)^2,\\
&&\bcm\, E_3\le \frac{_{8}}{^3}\zeta^2(3/2)\pi^3\sigma m^{3/2}N^{-1/2}\log N d(N)\,e^{-N/(2\pi m\sigma\log N)}, \\
&&\bcm\, E^{(3)}\le 16\pi^3m\sum_{c>0, N|c}min[\frac{_2}{^\pi}\phi(q)c^{-1}\log c,\frac{_1}{^6}\sigma N\log N m^{1/2}c^{-3/2}d(c)].
\end{eqnarray*}
\end{thm}
As we said,  we also need the following explicit bound
\begin{lem}(Ellenberg)\label{ellen2} Let $p$ be a prime, $m$ a positive integer, $\chi$ a quadratic character
of conductor $q$ prime to $p$. Then
$$
(a_{m},L_{\chi})_p\le 2\sqrt{3}m^{1/2}\tau(m)(1-e^{-2\pi/q\sqrt{p}})^{-1}
(4\pi+16\zeta(3/2)^2\pi^2p^{-3/2}).
$$
\end{lem}
{\bf Remark:} In \cite{e2} the author only proves a version of this lemma valid only for $m$ a multiple of $p$. However, as he himself mention in \cite{e3} the result
is valid as stated in Lemma \ref{ellen2}. 

\medskip

With these two results on hand we can prove the following.
\begin{lem}\label{analitischen} We have
$
|(a_1,L_\chi)_{p^2}^{p-\text{new}}|\ne 0
$
whenever 
\begin{itemize}
 \item $p> 349$ and  $\chi$ is a Dirichlet character of order at most  $8$.
\item $p> 131$ and  $\chi$ is a Dirichlet character of order at most  $3$.
\end{itemize}

\end{lem}
The proof of the lemma is a direct application of  Theorem
\ref{ellen} and Lemma \ref{ellen2} to our case. We sketch the proof for the first part of the lemma and, hence, we suppose $q\le 8$. The proof for the second part is  exactly the same  but just changing the value $q$ of the conductor. In order to use Theorem \ref{ellen} we will choose $\sigma=q^2/2\pi$. With those values let $p\ge 457$ be a
given prime and $x=\frac{{32}}{\pi} N\log N$ where $N$ will be
either $p$ or $p^2$. Hence we get, for $m=1$, $N=p^2$,
\begin{eqnarray*}
&&\left|(a_1,B(x))_p\right|\le 5760 \,e^{2\pi}(400/399)^{3}p^{-3},\\
&&E_1\le\frac{_{1024}}{^3}\pi^2\log p \,e^{-p^2/(128\log p)},\\
&&E_2\le \frac{_{65536}}{^9}\pi^3p^{-3}(\log p)^2\zeta(7/2)^2,\\
 &&E_3\le \frac{_{1024}}{^3}\zeta(3/2)^2\pi^2p^{-1}\,e^{-p^2/128\log p}\log p, \\
 &&E^{(3)}\le 1280\pi^2(\log p)^2p^{-2} +\frac{_{1536}}{^3}\pi^2\log p\, p^{-1}
 \left(\zeta^2(3/2)-\sum_{k=1}^p\frac{\tau(k)}{k^{3/2}}\right).
\end{eqnarray*}
In order to get the bound for $E^{(3)}$  we have split the sum in Theorem \ref{ellen} depending on whether $c>p^2$ or not.
We just have to note that all the functions appearing are decreasing in $p$ to get, with the help of maple,
$$
|(a_1,L_\chi)_{p^2}-4\pi|\le 9\text{.} 16\,,
$$
by  substituting $p=353$. On the other hand, the same argument, but now using Lemma \ref{ellen2} instead gives us
\begin{eqnarray*}
&&\frac{1}{p^2-1}|(a_p,L_\chi)_p|\le 0\text{.}33\,\qquad\text{and}\\
&&\frac{p}{p^2-1}|(a_1,L_\chi)_p|\le 3.06.
\end{eqnarray*}
for all $p\ge 353$. The result follows by (\ref{inner}). 

\medskip

Substituting $q$ by $3$ instead of $8$ we get the bounds
\begin{eqnarray*}
&&|(a_1,L_\chi)_{p^2}-4\pi|\le 9\text{.} 6\,,\\
&&\frac{1}{p^2-1}|(a_p,L_\chi)_p|\le 0\text{.}35\,\\
&&\frac{p}{p^2-1}|(a_1,L_\chi)_p|\le 2.05,
\end{eqnarray*}
for any $p\ge 137$ which gives the second part of the lemma.

\end{document}